\documentclass[11pt, reqno]{amsart}%
\usepackage{amsmath}%
\usepackage{amsfonts}%
\usepackage{amssymb}%
\usepackage{graphicx}

\newtheorem{theorem}{Theorem}
\theoremstyle{plain}

\newtheorem{corollary}{Corollary}

\newtheorem{definition}{Definition}
\newtheorem{example}{Example}

\newtheorem{proposition}{Proposition}
\newtheorem{remark}{Remark}

\numberwithin{equation}{section}

\begin{document}
\title[Pseudo-stopping times]{A definition and some characteristic properties of pseudo-stopping times}
\author{Ashkan Nikeghbali}
\address{Laboratoire de Probabilit\'{e}s et Mod\`{e}les Al\'{e}atoires\\
Universit\'{e} Pierre et Marie Curie, et CNRS UMR 7599\\
175, rue du Chevaleret F-75013 Paris, France \\
and C.R.E.S.T\\
( Centre de Recherche en Economie et Statistique Th\'{e}orique)}
\email{nikeghba@ccr.jussieu.fr}
\author{Marc Yor}
\address{Laboratoire de Probabilit\'{e}s et Mod\`{e}les Al\'{e}atoires\\
Universit\'{e} Pierre et Marie Curie, et CNRS UMR 7599\\
175, rue du Chevaleret F-75013 Paris, France}
\date{November 19, 2004}
\begin{abstract}
Recently, D. Williams \cite{williams} gave an explicit example of a random
time $\rho$ associated with Brownian motion such that $\rho$\ is not a
stopping time but $\mathbb{E}M_{\rho}=\mathbb{E}M_{0}$ for every bounded
martingale $M$. The aim of this paper is to give some characterizations for
such random times, which we call pseudo-stopping times, and to construct
further examples, using techniques of progressive enlargements of filtrations.

\end{abstract}
\maketitle


\section{Introduction}

Let $\left(  \Omega,\mathcal{F},\left(  \mathcal{F}_{t}\right)  _{t\geq
0},\mathbb{P}\right)  $ be a filtered probability space, and $\rho:$ $\left(
\Omega,\mathcal{F}\right)  \rightarrow\left(  \mathbb{R}_{+},\mathcal{B}%
\left(  \mathbb{R}_{+}\right)  \right)  $ be a random time. We recall that the
space $\mathcal{H}^{1}$ is the Banach space of (c\`{a}dl\`{a}g) $\left(
\mathcal{F}_{t}\right)  $-martingales $\left(  M_{t}\right)  $ such that
\[
\left\Vert M\right\Vert _{\mathcal{H}^{1}}=\mathbb{E}\left[  \sup_{t\geq
0}\left\vert M_{t}\right\vert \right]  <\infty.
\]

\begin{definition}
We say that $\rho$ is a $\left(  \mathcal{F}_{t}\right)  $ pseudo-stopping
time if for every$\ \left(  \mathcal{F}_{t}\right)  $-martingale $\left(
M_{t}\right)  $ in $\mathcal{H}^{1}$, we have%
\begin{equation}
\mathbb{E}M_{\rho}=\mathbb{E}M_{0}.\label{pta}%
\end{equation}
\end{definition}

\begin{remark}
It is equivalent to assume that (\ref{pta}) holds for bounded martingales,
since these are dense in $\mathcal{H}^{1}$.
\end{remark}

\bigskip

We indicate immediately that a class of pseudo-stopping times with respect to
a filtration $\left(  \mathcal{F}_{t}\right)  $\ which are not in general
$\left(  \mathcal{F}_{t}\right)  $\ stopping times may be obtained by
considering stopping times with respect to a larger filtration $\left(
\mathcal{G}_{t}\right)  $\ such that $\left(  \mathcal{F}_{t}\right)  $ is
immersed in $\left(  \mathcal{G}_{t}\right)  $, i.e: every $\left(
\mathcal{F}_{t}\right)  $\ martingale is a $\left(  \mathcal{G}_{t}\right)
$\ martingale. This situation is described in (\cite{bremaudyor}) and refered
to there as the $\left(  H\right)  $\ hypothesis. We shall discuss this
situation in more details in Section 3. For now, we give an example. Let
$B_{t}=\left(  B_{t}^{1},\ldots,B_{t}^{d}\right)  $\ be a $d$-dimensional
Brownian motion, and $R_{t}=\left\vert B_{t}\right\vert $, $t\geq0,$ its
radial part; it is well known that
\[
\left(  \mathcal{R}_{t}\equiv\sigma\left\{  R_{s},\text{ }s\leq t\right\}
,\text{ }t\geq0\right)  ,
\]
the natural filtration of $R$,\ is immersed in $\left(  \mathcal{B}_{t}%
\equiv\sigma\left\{  B_{s},\text{ }s\leq t\right\}  ,\text{ }t\geq0\right)  $,
the natural filtration of $B$. Thus an example of $\left(  \mathcal{R}%
_{t}\right)  $\ pseudo-stopping time is:%
\[
T_{a}^{\left(  1\right)  }=\inf\left\{  t,\text{ }B_{t}^{1}>a\right\}  .
\]

Recently, D. Williams \cite{williams} showed that with respect to the
filtration $\left(  \mathcal{F}_{t}\right)  $\ generated by a one dimensional
Brownian motion $\left(  B_{t}\right)  _{t\geq0}$, there exist pseudo-stopping
times $\rho$\ which are not $\left(  \mathcal{F}_{t}\right)  $\ stopping
times. D. Williams' example is the following: let
\[
T_{1}=\inf\left\{  t:\text{ }B_{t}=1\right\}  ,\text{ }\sigma=\sup\left\{
t<T_{1}:\text{ }B_{t}=0\right\}  ;
\]
then%
\[
\rho=\sup\left\{  s<\sigma:\text{ }B_{s}=S_{s}\right\}  ,\text{ where }%
S_{s}=\sup_{u\leq s}B_{u}%
\]
is a $\left(  \mathcal{F}_{t}\right)  $\ pseudo-stopping time. This paper has
two main aims:

\begin{itemize}
\item to understand better the nature of pseudo-stopping times;

\item to construct further examples of pseudo-stopping times;\bigskip
\end{itemize}

In Section 2, with the help of the theory of progressive enlargements of
filtrations, we give some equivalent properties for $\rho$\ to be a
pseudo-stopping time. We also comment there on the difference between
(\ref{pta}) and the property
\begin{equation}
\mathbb{E}\left[  M_{\infty}\mid\mathcal{F}_{\rho}\right]  =M_{\rho
}\label{knimais}%
\end{equation}
for every uniformly integrable $\left(  \mathcal{F}_{t}\right)  $-martingale
$\left(  M_{t}\right)  $, which was shown by Knight and Maisonneuve
\cite{Knight-Maisonneuve} to be equivalent to $\rho$\ being a $\left(
\mathcal{F}_{t}\right)  $-stopping time.\bigskip

In Section 3, we give some other examples of pseudo-stopping times. We
associate with the end $L$\ of a given $\left(  \mathcal{F}_{t}\right)
$\ predictable set $\Gamma$, i.e
\[
L=\sup\left\{  t:\left(  t,\omega\right)  \in\Gamma\right\}  ,
\]
a pseudo-stopping time $\rho<L$ in a manner which generalizes D. Williams'
example. We also link the pseudo-stopping times with randomized stopping times.\bigskip

In Section 4, we give a discrete time analogue of the Williams random time
$\rho$. This approach is based on the analogue of Williams' path decomposition
obtained by Le Gall for the standard random walk \cite{Le Gall}.

\section{Some characteristic properties of pseudo-stopping times}

\subsection{Basic facts about progressive enlargements}

We recall here some basic results about the progressive enlargement of a
filtration $\left(  \mathcal{F}_{t}\right)  $\ by a random time $\rho$. All
these results may be found in \cite{jeulinyor}, \cite{jeulin}, \cite{zurich},
\cite{delmaismey} or \cite{protter}.\bigskip

We enlarge the initial filtration $\left(  \mathcal{F}_{t}\right)  $\ with the
process $\left(  \rho\wedge t\right)  _{t\geq0}$, so that the new enlarged
filtration $\left(  \mathcal{F}_{t}^{\rho}\right)  _{t\geq0}$\ is the smallest
filtration containing $\left(  \mathcal{F}_{t}\right)  $\ and making $\rho$\ a
stopping time. A few processes play a crucial role in our discussion:

\begin{itemize}
\item the $\left(  \mathcal{F}_{t}\right)  $-supermartingale
\begin{equation}
Z_{t}^{\rho}=\mathbb{P}\left[  \rho>t\mid\mathcal{F}_{t}\right]
\label{surmart}%
\end{equation}
chosen to be c\`{a}dl\`{a}g, associated to $\rho$\ by Az\'{e}ma (see
\cite{jeulin} for detailed references);

\item the $\left(  \mathcal{F}_{t}\right)  $-dual optional and predictable
projections of the process $1_{\left\{  \rho\leq t\right\}  }$, denoted
respectively by $A_{t}^{\rho}$ and $a_{t}^{\rho}$;

\item the c\`{a}dl\`{a}g martingale
\[
\mu_{t}^{\rho}=\mathbb{E}\left[  A_{\infty}^{\rho}\mid\mathcal{F}_{t}\right]
=A_{t}^{\rho}+Z_{t}^{\rho}%
\]
which is in BMO($\mathcal{F}_{t}$) (see \cite{delmaismey} or \cite{zurich}).
We recall that the space of BMO martingales (see \cite{dellachmeyer} for more
details and references) is the Banach space of (c\`{a}dl\`{a}g) square
integrable $\left(  \mathcal{F}_{t}\right)  $-martingales $\left(
Y_{t}\right)  $ which satisfy%
\[
\left\Vert Y\right\Vert _{BMO}^{2}=\mbox{essup} _{T}\mathbb{E}\left[  \left(
Y_{\infty}-Y_{T-}\right)  ^{2}\mid\mathcal{F}_{T}\right]  <\infty
\]
where $T$\ ranges over all $\left(  \mathcal{F}_{t}\right)  $-stopping times.
\end{itemize}

We also consider the Doob-Meyer decomposition of (\ref{surmart}):%
\[
Z_{t}^{\rho}=m_{t}^{\rho}-a_{t}^{\rho}%
\]
If $\rho$\ avoids any $\left(  \mathcal{F}_{t}\right)  $-stopping time, that
is to say $P\left[  \rho=T>0\right]  =0$ for any stopping time $T$, then
$A_{t}^{\rho}=a_{t}^{\rho}$ is continuous.\bigskip

Finally, we recall that every $\left(  \mathcal{F}_{t}\right)  $-local
martingale $\left(  M_{t}\right)  $, stopped at $\rho$, is a $\left(
\mathcal{F}_{t}^{\rho}\right)  $\ semimartingale, with canonical
decomposition:
\begin{equation}
M_{t\wedge\rho}=\widetilde{M}_{t}+\int_{0}^{t\wedge\rho}\frac{d<M,\mu^{\rho
}>_{s}}{Z_{s-}^{\rho}}\label{decocanonique}%
\end{equation}
where $\left(  \widetilde{M}_{t}\right)  $\ is an $\left(  \mathcal{F}%
_{t}^{\rho}\right)  $-local martingale.

\begin{remark}
We also recall that in a filtration $\left(  \mathcal{F}_{t}\right)  $\ where
all martingales are continuous, $A_{t}^{\rho}=a_{t}^{\rho}$ since optional
processes are predictable (see \cite{revuzyor}, chapter IV).
\end{remark}

\subsection{A characterization of pseudo-stopping times}

We now discuss some characteristic properties of pseudo-stopping times. We
assume throughout that $\mathbb{P}\left[  \rho=\infty\right]  =0$.\

\begin{theorem}
\label{mainthm}The following four properties are equivalent:

\begin{enumerate}
\item $\rho$\ is a $\left(  \mathcal{F}_{t}\right)  $\ pseudo-stopping time,
i.e (\ref{pta}) is satisfied;

\item $\mu_{t}^{\rho}\equiv1$, $a.s$

\item $A_{\infty}^{\rho}\equiv1$, $a.s$

\item every $\left(  \mathcal{F}_{t}\right)  $\ local martingale $\left(
M_{t}\right)  $ satisfies
\[
\left(  M_{t\wedge\rho}\right)  _{t\geq0}\text{ }is\text{ }a\text{
}local\text{ }\left(  \mathcal{F}_{t}^{\rho}\right)  \text{ }martingale.
\]

If, furthermore, all $\left(  \mathcal{F}_{t}\right)  $ martingales are
continuous, then each of the preceding properties is equivalent to

\item
\[
\left(  Z_{t}^{\rho}\right)  _{t\geq0}\text{ is a decreasing }\left(
\mathcal{F}_{t}\right)  \ \text{predictable process}%
\]
\end{enumerate}
\end{theorem}

\begin{proof}
$\left(  1\right)  \Longrightarrow\left(  2\right)  $ For every square
integrable $\left(  \mathcal{F}_{t}\right)  $\ martingale $\left(
M_{t}\right)  $, we have
\[
\mathbb{E}\left[  M_{\rho}\right]  =\mathbb{E}\left[  \int_{0}^{\infty}%
M_{s}dA_{s}^{\rho}\right]  =\mathbb{E}\left[  M_{\infty}A_{\infty}^{\rho
}\right]  =\mathbb{E}\left[  M_{\infty}\mu_{\infty}^{\rho}\right]  .
\]
Since $\mathbb{E}M_{\rho}=\mathbb{E}M_{0}=\mathbb{E}M_{\infty}$, we have
\[
\mathbb{E}\left[  M_{\infty}\right]  =\mathbb{E}\left[  M_{\infty}A_{\infty
}^{\rho}\right]  =\mathbb{E}\left[  M_{\infty}\mu_{\infty}^{\rho}\right]  .
\]
Consequently, $\mu_{\infty}^{\rho}\equiv1$, $a.s$,$\ $hence $\mu_{t}^{\rho
}\equiv1$, $a.s$ which is equivalent to: $A_{\infty}^{\rho}\equiv1$, $a.s$.
Hence, $2.$ and $3.$ are equivalent.

$\left(  2\right)  \Longrightarrow\left(  4\right)  .$ This is a consequence
of the decomposition formula (\ref{decocanonique}).

$\left(  4\right)  \Longrightarrow\left(  1\right)  .$ It suffices
to consider any $\mathcal{H}^{1}$-martingale $\left(  M_{t}\right)
$, which, assuming $4.$, satisfies: $\left(  M_{t\wedge\rho}\right)
_{t\geq0}$ is a martingale in the enlarged filtration, for which
$\rho$\ is a stopping time. Then as a consequence of the optional
stopping theorem applied in $\left(
\mathcal{F}_{t}^{\rho}\right)  $ at time $\rho$, we get%
\[
\mathbb{E}\left[  M_{\rho}\right]  =\mathbb{E}\left[  M_{0}\right]  ,
\]
hence $\rho\ $is a pseudo-stopping time.

Finally, in the case where all $\left(  \mathcal{F}_{t}\right)  $\ martingales
are continuous, we show:

$\left.  a\right)  $ $\left(  2\right)  \Rightarrow\left(  5\right)  $ If
$\rho$\ is a pseudo-stopping time, then $Z_{t}^{\rho}$ decomposes as%
\[
Z_{t}^{\rho}=1-A_{t}^{\rho}.
\]
As all $\left(  \mathcal{F}_{t}\right)  $\ martingales are continuous,
optional processes are in fact predictable, and so $\left(  Z_{t}^{\rho
}\right)  $\ is a predictable decreasing process.

$\left.  b\right)  $ $\left(  5\right)  \Rightarrow\left(  2\right)  $
Conversely, if $\left(  Z_{t}^{\rho}\right)  $\ is a predictable decreasing
process, then from the unicity in the Doob-Meyer decomposition, the martingale
part $\mu_{t}^{\rho}$ is constant, i.e. $\mu_{t}^{\rho}\equiv1$, $a.s$. Thus,
$\left(  2\right)  $ is satisfied.
\end{proof}

In the next proposition, we deal with uniformly integrable martingales
$\left(  M_{t}\right)  $ instead of martingales in $\mathcal{H}^{1}$ (or
$\mathcal{H}^{2},\ldots$).

\begin{proposition}
The following properties are equivalent:

\begin{enumerate}
\item $\rho$\ is a $\left(  \mathcal{F}_{t}\right)  $\ pseudo-stopping time;

\item for every uniformly integrable martingale%
\[
\mathbb{E}\left[  \left\vert M_{\rho}\right\vert \right]  \leq\mathbb{E}%
\left[  \left\vert M_{\infty}\right\vert \right]  .
\]
\end{enumerate}
\end{proposition}

\begin{remark}
In fact, we shall further show in the next proof, that for $\rho$\ a
pseudo-stopping time and for $\left(  M_{t}\right)  $\ a uniformly integrable
martingale:%
\[
\mathbb{E}\left[  \left\vert M_{\rho}\right\vert \right]  <\infty,\text{ and
}\mathbb{E}\left[  M_{\rho}\right]  =\mathbb{E}\left[  M_{\infty}\right]  .
\]
\end{remark}

\begin{proof}
$\left(  1\right)  \Rightarrow\left(  2\right)  $ If $\left(  M_{t}\right)  $
is uniformly integrable, it may be decomposed as:%
\begin{equation}
M_{t}=M_{t}^{\left(  +\right)  }-M_{t}^{\left(  -\right)  }\label{posnegpart}%
\end{equation}
where
\[
M_{t}^{\left(  +\right)  }=\mathbb{E}\left[  M_{\infty}^{+}\mid\mathcal{F}%
_{t}\right]  \text{ and }M_{t}^{\left(  -\right)  }=\mathbb{E}\left[
M_{\infty}^{-}\mid\mathcal{F}_{t}\right]  .
\]
(Note that $M_{\infty}^{\pm}$\ indicates the positive and negative parts of
$M_{\infty}$, whereas $\left(  M_{t}^{\left(  \pm\right)  }\right)  $\ are the
martingales with terminal values $M_{\infty}^{\pm}$). Thus to prove $\left(
2\right)  $ it suffices to prove%
\[
\mathbb{E}\left[  M_{\rho}\right]  =\mathbb{E}\left[  M_{\infty}\right]  ,
\]
under the further assumption that $M\geq0$. In this latter case, we have
$M_{t}=\mathbb{E}\left[  M_{\infty}\mid\mathcal{F}_{t}\right]  $, with
$M_{\infty}\geq0$. Now let
\[
M_{t}^{\left(  n\right)  }=\mathbb{E}\left[  \left(  M_{\infty}\wedge
n\right)  \mid\mathcal{F}_{t}\right]  .
\]
$\left(  M_{t}^{\left(  n\right)  }\right)  $ is a bounded martingale, hence
we have%
\[
\mathbb{E}\left[  M_{\infty}^{\left(  n\right)  }\right]  =\mathbb{E}\left[
M_{\rho}^{\left(  n\right)  }\right]  .
\]
We also have
\[
\mathbb{P}\left[  \sup_{t\geq0}\left(  M_{t}-M_{t}^{\left(  n\right)
}\right)  >\varepsilon\right]  \leq\frac{1}{\varepsilon}\mathbb{E}\left[
M_{\infty}-M_{\infty}^{\left(  n\right)  }\right]  ,
\]
so that $\left(  M_{\rho}^{\left(  n\right)  }\right)  $ converges to $\left(
M_{\rho}\right)  $ in probability; but the sequence $\left(  M_{\rho}^{\left(
n\right)  }\right)  $\ is increasing, so it in fact converges almost surely.
Hence the monotone convergence theorem yields
\[
\mathbb{E}\left[  M_{\infty}\right]  =\mathbb{E}\left[  M_{\rho}\right]  .
\]
Finally, going back to (\ref{posnegpart}) in the general case, we obtain:%
\begin{align*}
\mathbb{E}\left[  \left\vert M_{\rho}\right\vert \right]   & \leq
\mathbb{E}\left[  M_{\rho}^{\left(  +\right)  }+M_{\rho}^{\left(  -\right)
}\right] \\
& =\mathbb{E}\left[  M_{\infty}^{+}+M_{\infty}^{-}\right] \\
& =\mathbb{E}\left[  \left\vert M_{\infty}\right\vert \right]  .
\end{align*}
Hence $\left(  2\right)  $ holds. Further, we may now write:%
\begin{align*}
\mathbb{E}\left[  M_{\rho}\right]   & =\mathbb{E}\left[  M_{\rho}^{\left(
+\right)  }-M_{\rho}^{\left(  -\right)  }\right] \\
& =\mathbb{E}\left[  M_{\infty}^{+}-M_{\infty}^{-}\right] \\
& =\mathbb{E}\left[  M_{\infty}\right]  .
\end{align*}

$\left(  2\right)  \Rightarrow\left(  1\right)  $ We need only apply
property $\left(  2\right)  $ to any martingale $\left(
M_{t}\right)  $ taking values
in $\left[  0,1\right]  $. Thus:%
\begin{align*}
\mathbb{E}\left[  M_{\rho}\right]   & \leq\mathbb{E}\left[  M_{\infty}\right]
\\
\mathbb{E}\left[  1-M_{\rho}\right]   & \leq\mathbb{E}\left[  1-M_{\infty
}\right]  .
\end{align*}
But, since the sums on both sides add up to $1$, we must have:%
\[
\mathbb{E}\left[  M_{\rho}\right]  =\mathbb{E}\left[  M_{\infty}\right]
\]
Hence, $\rho$ is a $\left(  \mathcal{F}_{t}\right)  $\ pseudo-stopping time.
\end{proof}

As an application of the theorem, we can check that in D. Williams' example,
his time $\rho$ associated with a Brownian motion is a pseudo-stopping time.
Indeed, the dual predictable (=optional) projection $A_{t}^{\rho}$\ of
$1_{\left\{  \rho\leq t\right\}  }$ is max$_{s\leq t\wedge T_{1}}B_{s}$
(\cite{williams}, \cite{rogerswilliams}) and $A_{\infty}^{\rho}\equiv1$.

\subsection{Around the result of Knight and Maisonneuve}

We now comment on the statement of the fourth property in Theorem
\ref{mainthm}.

For the properties of the different sigma fields $\mathcal{F}_{\rho}$,
$\mathcal{F}_{\rho+}$, $\mathcal{F}_{\rho-}$, associated with a general random
time $\rho$, the reader can consult \cite{rogerswilliams} or \cite{zurich}.
Here, we just recall the definitions:

\begin{definition}
Three classical $\sigma$-fields associated with a filtration $\left(
\mathcal{F}_{t}\right)  \ $and any random time $\rho$ are:
\[
\left\{
\begin{array}
[c]{lcl}%
\mathcal{F}_{\rho+} & = & \sigma\left\{  z_{\rho},\text{ }\left(
z_{t}\right)  \text{ any }\left(  \mathcal{F}_{t}\right)
\ \text{progressively measurable process}\right\}  ;\\
\mathcal{F}_{\rho} & = & \sigma\left\{  z_{\rho},\text{ }\left(  z_{t}\right)
\text{ any }\left(  \mathcal{F}_{t}\right)  \ \text{optional process}\right\}
;\\
\mathcal{F}_{\rho-} & = & \sigma\left\{  z_{\rho},\text{ }\left(
z_{t}\right)  \text{ any }\left(  \mathcal{F}_{t}\right)  \ \text{predictable
process}\right\}  ;
\end{array}
\right.
\]
\end{definition}

The result of Knight and Maisonneuve which was recalled in the introduction
may be stated as follows:

\begin{theorem}
If for all uniformly integrable $\left(  \mathcal{F}_{t}\right)  $-martingales
$\left(  M_{t}\right)  $, one has
\[
\mathbb{E}\left[  M_{\infty}\mid\mathcal{F}_{\rho}\right]  =M_{\rho},\qquad
on\text{ }\left\{  \rho<\infty\right\}  ,
\]
then $\rho$ is a $\left(  \mathcal{F}_{t}\right)  $-stopping time (the
converse is Doob's optional stopping theorem).
\end{theorem}

Refining slightly the argument in \cite{Knight-Maisonneuve}, we obtain the following:

\begin{theorem}
If for all bounded $\left(  \mathcal{F}_{t}\right)  $-martingales $\left(
M_{t}\right)  $, one has
\[
\mathbb{E}\left[  M_{\infty}\mid\mathcal{\sigma}\left\{  M_{\rho}%
,\rho\right\}  \right]  =M_{\rho},\qquad on\text{ }\left\{  \rho
<\infty\right\}  ,
\]
then $\rho$ is a $\left(  \mathcal{F}_{t}\right)  $-stopping time.
\end{theorem}

\begin{proof}
For $t\geq0$\ we have
\[
\mathbb{E}\left[  M_{\infty}\mathbf{1}_{\left(  \rho\leq t\right)  }\right]
=\mathbb{E}\left[  M_{\rho}1_{\left(  \rho\leq t\right)  }\right]
=\mathbb{E}\left[  \int_{0}^{t}M_{s}dA_{s}^{\rho}\right]  =\mathbb{E}\left[
M_{\infty}A_{t}^{\rho}\right]  .
\]
Comparing the two extreme terms, we get
\[
\mathbf{1}_{\left(  \rho\leq t\right)  }=A_{t}^{\rho},
\]
i.e $\rho$\ is a $\left(  \mathcal{F}_{t}\right)  $-stopping time.\bigskip
\end{proof}

An interesting open question in view of what has been proved for
pseudo-stopping times is whether $\mathbb{E}\left[  M_{\infty}\mid M_{\rho
}\right]  =M_{\rho},$ $on$ $\left\{  \rho<\infty\right\}  $\ is equivalent to
$\rho$\ being a stopping time.\bigskip

To illustrate the result of Knight and Maisonneuve, we show explicitly how, in
the framework of D. Williams' example, $M_{\rho}$\ and $\mathbb{E}\left[
M_{\infty}\mid\mathcal{F}_{\rho}\right]  $\ differ, for
\[
M_{t}=\exp\left(  \lambda B_{t\wedge T_{1}}-\frac{\lambda^{2}}{2}\left(
t\wedge T_{1}\right)  \right)  ,\qquad\lambda>0.
\]
We write%
\begin{align*}
M_{\infty}  & =\exp\left(  \lambda-\frac{\lambda^{2}}{2}T_{1}\right) \\
& =\exp\left(  \lambda\right)  \exp\left(  -\frac{\lambda^{2}}{2}\left(
\rho+\left(  \sigma-\rho\right)  +\left(  T_{1}-\sigma\right)  \right)
\right)
\end{align*}
and we compute:%
\begin{align}
\mathbb{E}\left[  M_{\infty}\mid\mathcal{F}_{\rho}\right]   & =\label{trois}\\
& \exp\left(  \lambda-\frac{\lambda^{2}}{2}\rho\right)  \mathbb{E}\left[
\exp\left(  -\frac{\lambda^{2}}{2}\left(  \sigma-\rho\right)  \right)
\mid\mathcal{F}_{\rho}\right]  \mathbb{E}\left[  \exp\left(  -\frac{\lambda
^{2}}{2}\left(  T_{1}-\sigma\right)  \right)  \right]  ,
\end{align}
since $\left(  T_{1}-\sigma\right)  $ is independent from $\mathcal{F}%
_{\sigma}$, (and consequently from $\mathcal{F}_{\rho}$, since $\mathcal{F}%
_{\rho}\subset\mathcal{F}_{\sigma}$).

We now recall D. Williams' path decomposition results for $\left(
B_{u}\right)  _{u\leq T_{1}}$ on the intervals $\left(  0,\rho\right)  $,
$\left(  \rho,\sigma\right)  $, $\left(  \sigma,T_{1}\right)  $:

\begin{itemize}
\item $\left(  B_{\sigma+u}\right)  _{u\leq T_{1}-\sigma}$ is a BES$\left(
3\right)  $ process, independent of $\mathcal{F}_{\sigma}$; hence we have
\[
\mathbb{E}\left[  \exp\left(  -\frac{\lambda^{2}}{2}\left(  T_{1}%
-\sigma\right)  \right)  \right]  =\frac{\lambda}{\sinh\left(  \lambda\right)
}.
\]

\item $S_{\rho}$, where $S_{s}=\sup_{u\leq s}B_{u}$, is uniformly distributed
on $\left(  0,1\right)  $;

\item Conditionally on $S_{\rho}=h$, the processes $\left(  B_{u}\right)
_{u\leq\rho}$ and $\left(  B_{\sigma-u}\right)  _{u\leq\sigma-\rho}$ are two
independent Brownian motions considered up to their first hitting time of $h$.
Consequently, we have:%
\[
\mathbb{E}\left[  \exp\left(  -\frac{\lambda^{2}}{2}\left(  \sigma
-\rho\right)  \right)  \mid\mathcal{F}_{\rho}\right]  =\exp\left(  -\lambda
S_{\rho}\right)  .
\]
Plugging these informations in (\ref{trois}), we obtain:%
\[
\mathbb{E}\left[  M_{\infty}\mid\mathcal{F}_{\rho}\right]  =\exp\left(
\lambda\left(  1-B_{\rho}\right)  -\frac{\lambda^{2}}{2}\rho\right)  \left(
\frac{\lambda}{\sinh\left(  \lambda\right)  }\right)  ,
\]
whilst%
\begin{equation}
M_{\rho}=\exp\left(  \lambda B_{\rho}-\frac{\lambda^{2}}{2}\rho\right)
\end{equation}
and these two quantities are obviously different.
\end{itemize}

\subsection{Further properties of pseudo-stopping times}

Besides the assumption that $\rho$\ is a $\left(  \mathcal{F}_{t}\right)  $
pseudo-stopping time, we also make the hypothesis that $\rho$ avoids all
$\left(  \mathcal{F}_{t}\right)  $-stopping times. We saw that in this case
\[
a_{t}^{\rho}=A_{t}^{\rho}=1-Z_{t}^{\rho}%
\]
is continuous.

For simplicity, we shall write $\left(  Z_{u}\right)  $\ instead of $\left(
Z_{u}^{\rho}\right)  $.

\begin{proposition}
\label{terminaval}Under the previous hypotheses, for all uniformly integrable
$\left(  \mathcal{F}_{t}\right)  $\ martingales $\left(  M_{t}\right)  $, and
all bounded Borel measurable functions $f$, one has:%
\[
\mathbb{E}\left[  M_{\rho}f\left(  Z_{\rho}\right)  \right]  =\mathbb{E}%
\left[  M_{0}\right]  \int_{0}^{1}f\left(  x\right)  dx=\mathbb{E}\left[
M_{\rho}\right]  \int_{0}^{1}f\left(  x\right)  dx.
\]
\end{proposition}

\begin{remark}
On the other hand it is not true that
\begin{equation}
\mathbb{E}\left[  M_{\infty}f\left(  Z_{\rho}\right)  \right]  =\mathbb{E}%
\left[  M_{\rho}f\left(  Z_{\rho}\right)  \right]  ,\label{rhs}%
\end{equation}
for every bounded Borel function $f$. Indeed, from Proposition
\ref{terminaval}, the right hand side of (\ref{rhs}) is equal to:%
\[
\mathbb{E}\left[  M_{\infty}\int_{0}^{1}f\left(  x\right)  dx\right]  .
\]
Thus, our hypothesis (\ref{rhs}) would imply the absurd equality between
$f\left(  Z_{\rho}\right)  $\ and $\int_{0}^{1}f\left(  x\right)  dx$.
\end{remark}

\begin{proof}
(of Proposition \ref{terminaval}) Under our assumptions, we have
\begin{align*}
\mathbb{E}\left[  M_{\rho}f\left(  Z_{\rho}\right)  \right]   & =\mathbb{E}%
\left[  \int_{0}^{\infty}M_{u}f\left(  Z_{u}\right)  dA_{u}^{\rho}\right] \\
& =\mathbb{E}\left[  \int_{0}^{\infty}M_{u}f\left(  1-A_{u}^{\rho}\right)
dA_{u}^{\rho}\right] \\
& =\mathbb{E}\left[  M_{\infty}\int_{0}^{\infty}f\left(  1-A_{u}^{\rho
}\right)  dA_{u}^{\rho}\right] \\
& =\mathbb{E}\left[  M_{\infty}\int_{0}^{1}f\left(  1-x\right)  dx\right] \\
& =\mathbb{E}\left[  M_{\infty}\int_{0}^{1}f\left(  x\right)  dx\right]  .
\end{align*}
\end{proof}

Taking $M_{t}\equiv1$, we find that $\left(  Z_{\rho}\right)  $\ is uniformly
distributed on $\left(  0,1\right)  $, which is already known
(\cite{jeulinyor}, \cite{zurich}) since (recalling that $Z_{u}$\ is
decreasing)%
\[
Z_{\rho}=\inf_{u\leq\rho}Z_{u}.
\]
In fact we have a stronger result: under all changes of probability on
$\mathcal{F}_{\rho}$, of the form%
\[
d\mathbb{Q}=M_{\rho}d\mathbb{P}%
\]
where $\left(  M_{t}\right)  $\ is a positive uniformly integrable $\left(
\mathcal{F}_{t}\right)  $-martingale such that $\mathbb{E}\left[
M_{0}\right]  =1$, the law of $Z_{\rho}$\ (is unchanged and) is uniform.

\begin{corollary}
Under the assumptions of Proposition \ref{terminaval}, we have
\[
\mathbb{E}\left[  M_{\rho}\mid Z_{\rho}\right]  =\mathbb{E}\left[  M_{\rho
}\right]  =\mathbb{E}\left[  M_{0}\right]
\]
On the other hand, the quantity $\mathbb{E}\left[  M_{\infty}\mid Z_{\rho
}\right]  $\ is not easy to evaluate, as is seen with D. Williams' example,
and is different from $\mathbb{E}\left[  M_{\rho}\mid Z_{\rho}\right]  $.
Indeed, in this framework and with the already used notations:%
\[
\mathbb{E}\left[  M_{\infty}\mid Z_{\rho}\right]  =\exp\left(  \lambda\right)
\mathbb{E}\left[  \exp\left(  -\frac{\lambda^{2}}{2}T_{1}\right)  \mid
B_{\rho}\right]  .
\]
Decomposing again $T_{1}$ as $T_{1}=\rho+\left(  \sigma-\rho\right)  +\left(
T_{1}-\sigma\right)  $, and using D. Williams\textquotedblright\ path
decomposition, we obtain:%
\begin{align*}
\mathbb{E}\left[  M_{\infty}\mid Z_{\rho}\right]   & =\exp\left(
\lambda\right)  \left(  \frac{\lambda}{\sinh\left(  \lambda\right)  }\right)
\exp\left(  -\lambda B_{\rho}\right)  \mathbb{E}\left[  \exp\left(
-\frac{\lambda^{2}}{2}\rho\right)  \mid B_{\rho}\right] \\
& =\left(  \frac{2\lambda}{1-\exp\left(  -2\lambda\right)  }\right)
\exp\left(  -2\lambda B_{\rho}\right)  .
\end{align*}
\end{corollary}

\begin{corollary}
The family $\left\{  M_{\rho};\text{ }M\text{ uniformly integrable }\left(
\mathcal{F}_{t}\right)  \text{-martingale}\right\}  $ is not dense in
$L^{1}\left(  \mathcal{F}_{\rho}\right)  $.
\end{corollary}

\begin{proof}
From Proposition 2, the variable $\left(  f\left(  Z_{\rho}\right)  -\int
_{0}^{1}f\left(  x\right)  dx\right)  $ is orthogonal to $M_{\rho}$.
\end{proof}

This negative result led us to look for some representation of the generic
element of $L^{1}\left(  \mathcal{F}_{\rho}\right)  $\ in terms of $\left(
\mathcal{F}_{t}\right)  $-martingales taken at time $\rho$\ on one hand, and
the variable $Z_{\rho}$, on the other hand.

\begin{proposition}
$\left(  i\right)  .$ Let $K:$ $\left[  0,1\right]  \times\mathbb{R}_{+}%
\times\Omega\rightarrow\mathbb{R}_{+}$, be a $\mathcal{B}_{\left[  0,1\right]
}\otimes\mathcal{P}\left(  \mathcal{F}_{\bullet}\right)  $ measurable process,
where $\mathcal{P}\left(  \mathcal{F}_{\bullet}\right)  $ denotes the $\left(
\mathcal{F}_{t}\right)  $ predictable $\sigma$-field on $\mathbb{R}_{+}%
\times\Omega$. \ Then:%
\begin{equation}
\mathbb{E}\left[  K\left(  1-Z_{\rho},\rho\right)  \right]  =\mathbb{E}\left[
\int_{0}^{1}dyK\left(  y,\alpha_{y}\right)  \right] \label{masterformula}%
\end{equation}
where
\[
\alpha_{y}=\inf\left\{  u:\text{ }A_{u}^{\rho}>y\right\}  .
\]
$\qquad\qquad\qquad\qquad\quad\left(  ii\right)  .$ Let $\left(  H_{u}%
,u\geq0\right)  $ be a bounded predictable process. Define a measurable family
$\left(  M_{t}^{y}\right)  _{t\geq0}$ of martingales through their terminal
values:%
\[
M_{\infty}^{y}=H_{\alpha_{y}}.
\]
Then%
\[
H_{\rho}=M_{\rho}^{1-Z_{\rho}},\text{ }a.s.
\]
\end{proposition}

\begin{proof}
$\left(  i\right)  .$ This follows from the monotone class theorem, once we
have shown:%
\begin{equation}
\mathbb{E}\left[  f\left(  1-Z_{\rho}\right)  H_{\rho}\right]  =\mathbb{E}%
\left[  \int_{0}^{1}dyf\left(  y\right)  H_{\alpha_{y}}\right]
\label{formulerandom}%
\end{equation}
for every bounded predictable process $H$ and every Borel bounded function $f
$. But, this identity follows from the fact that: $1-Z_{\rho}=A_{\rho}$; and
so:%
\begin{align*}
\mathbb{E}\left[  f\left(  A_{\rho}\right)  H_{\rho}\right]   & =\mathbb{E}%
\left[  \int_{0}^{\infty}dA_{u}f\left(  A_{u}\right)  H_{u}\right] \\
& =\mathbb{E}\left[  \int_{0}^{1}dyf\left(  y\right)  H_{\alpha_{y}}\right]  .
\end{align*}
We shall prove the second statement by showing that for every bounded $\left(
k_{u}\right)  $ predictable process%
\[
\mathbb{E}\left[  k_{\rho}H_{\rho}\right]  =\mathbb{E}\left[  k_{\rho}M_{\rho
}^{1-Z_{\rho}}\right]  .
\]
From (\ref{masterformula}), we deduce:%
\begin{align*}
\mathbb{E}\left[  k_{\rho}M_{\rho}^{1-Z_{\rho}}\right]   & =\mathbb{E}\left[
\int_{0}^{1}dyM_{\alpha_{y}}^{y}k_{\alpha_{y}}\right] \\
& \overset{\left(  a\right)  }{=}\int_{0}^{1}dy\mathbb{E}\left[  M_{\infty
}^{y}k_{\alpha_{y}}\right] \\
& \overset{\left(  b\right)  }{=}\int_{0}^{1}dy\mathbb{E}\left[  H_{\alpha
_{y}}k_{\alpha_{y}}\right] \\
& \overset{\left(  c\right)  }{=}\mathbb{E}\left[  k_{\rho}H_{\rho}\right]  .
\end{align*}
((a) follows from the optional stopping theorem for $\left(  M_{t}^{y}\right)
$; (b) follows from the definition of $M_{\infty}^{y}$; (c) is another
consequence of (\ref{masterformula})). Comparing the extreme terms in the
above, we get
\[
H_{\rho}=M_{\rho}^{1-Z_{\rho}}.
\]
\end{proof}

\section{Some systematic constructions and some examples of pseudo-stopping times}

\subsection{First constructions}

Here we discuss some combinations of several pseudo-stopping times which yield
a pseudo-stopping time. Here is a first easy result:

\begin{proposition}
\label{construc}Let $\rho$\ be a $\left(  \mathcal{F}_{t}\right)
$-pseudo-stopping time and let $\tau$ be a $\left(  \mathcal{F}_{t}^{\rho
}\right)  $-stopping time. Then $\rho\wedge\tau$ is a $\left(  \mathcal{F}%
_{t}\right)  $ pseudo-stopping time.
\end{proposition}

\begin{proof}
Let $M$\ be any uniformly integrable $\left(  \mathcal{F}_{t}\right)
$-martingale. We know that $M_{t\wedge\rho}$ is a uniformly integrable
martingale in the enlarged filtration $\left(  \mathcal{F}_{t}^{\rho}\right)
$ and $\rho$\ is a stopping time in this filtration. If $\tau$\ is also a
$\left(  \mathcal{F}_{t}^{\rho}\right)  $-stopping time, then so is
$\rho\wedge\tau$. Hence $\mathbb{E}M_{\rho\wedge\tau}=\mathbb{E}M_{0}$.
\end{proof}

\begin{example}
Let $\rho$\ be as in D. Williams' example. Let $0<a<1$, and$\ T_{a}%
=\inf\left\{  t>0:\quad B_{t}=a\right\}  $. Then
\[
\rho_{a}=\rho\wedge T_{a},\qquad0<a<1,
\]
is an increasing family of pseudo-stopping times.
\end{example}

\begin{remark}
From the previous proposition, it is easy to see that for any uniformly
integrable $\left(  \mathcal{F}_{t}\right)  $-martingale, we have%
\[
\mathbb{E}\left[  M_{T\wedge\rho}\right]  =\mathbb{E}\left[  M_{0}\right]
\]
for any $\left(  \mathcal{F}_{t}\right)  $\ stopping time $T$.
\end{remark}

\begin{remark}
As a further comment about Proposition \ref{construc}, we remark that
pseudo-stopping times do not inherit all the "nice" properties of stopping
times. As an example, a pseudo-stopping time of a given filtration does not
remain in general a pseudo-stopping time in a larger filtration, whereas a
stopping time does. Indeed, let us keep the same notation as in section 2.3
and look at the pseudo-stopping time $\rho$\ in the larger filtration $\left(
\mathcal{F}_{t}^{\sigma}\right)  $. Using the computations we have already
done in section 2.3 and the projections formula (see \cite{delmaismey} p.186),
we get:%
\[
\mathbb{P}\left[  \rho>t\mid\mathcal{F}_{t}^{\sigma}\right]  =\frac{1-\max
_{s\leq t\wedge T_{1}}B_{s}}{1-B_{t\wedge T_{1}}^{+}}\mathbf{1}_{\left\{
\sigma>t\right\}  },
\]
which is not decreasing. In fact, any end of predictable set that avoids
stopping times is not a pseudo-stopping time, as we shall see in the next subsection.
\end{remark}

\subsection{A generalization of D. Williams' example}

To keep the discussion as simple as possible, we assume that we are working
with an original filtration $\left(  \mathcal{F}_{t}\right)  $\ such that:

\begin{itemize}
\item all $\left(  \mathcal{F}_{t}\right)  $-martingales are continuous (e.g:
$\left(  \mathcal{F}_{t}\right)  $ is the Brownian filtration).

\item Moreover, we consider $L$, the end of a $\left(  \mathcal{F}_{t}\right)
$\ predictable set, such that for every $\left(  \mathcal{F}_{t}\right)
$\ stopping time $T$, $\mathbb{P}\left[  L=T\right]  =0$.
\end{itemize}

Under these two conditions, the supermartingale $Z_{t}=P\left[  L>t\mid
\mathcal{F}_{t}\right]  $ associated with $L$\ is a.s. continuous, and
satisfies $Z_{L}=1$. Then we let,
\[
\rho=\sup\left\{  t<L:\quad Z_{t}=\inf_{u\leq L}Z_{u}\right\}  .
\]
The following holds:

\begin{proposition}
$\left(  i\right)  $ $I_{L}=\inf_{u\leq L}Z_{u}$ is uniformly distributed on
$\left[  0,1\right]  $; (see \cite{zurich})\newline

$\left(  ii\right)  $ The supermartingale $Z_{t}^{\rho}=P\left[  \rho
>t\mid\mathcal{F}_{t}\right]  $ associated with $\rho$\ is given by
\[
Z_{t}^{\rho}=\inf_{u\leq t}Z_{u}.
\]
As a consequence, $\rho$ is a $\left(  \mathcal{F}_{t}\right)  $
pseudo-stopping time.
\end{proposition}

\begin{proof}
$\left(  i\right)  $ Let
\[
T_{b}=\inf\left\{  t,\quad Z_{t}\leq b\right\}  ,\qquad0<b<1,
\]
then
\[
\mathbb{P}\left[  I_{L}\leq b\right]  =\mathbb{P}\left[  T_{b}<L\right]
=\mathbb{E}\left[  Z_{T_{b}}\right]  =b.
\]
$\left(  ii\right)  $ Note that for every $\left(  \mathcal{F}_{t}\right)  $
stopping time $T$, we have
\[
\left\{  T<\rho\right\}  =\left\{  T^{^{\prime}}<L\right\}
\]
where
\[
T^{^{\prime}}=\inf\left\{  t>T,\quad Z_{t}\leq\inf_{s\leq T}Z_{s}\right\}  .
\]
Consequently, we have%
\[
\mathbb{E}\left[  Z_{T}^{\rho}\right]  =\mathbb{P}\left[  T<\rho\right]
=\mathbb{P}\left[  T^{^{\prime}}<L\right]  =\mathbb{E}\left[  Z_{T^{^{\prime}%
}}\right]  =\mathbb{E}\left[  \inf_{u\leq T}Z_{u}\right]  .
\]
We deduce the desired result from the equality between the two extreme terms
for every $\left(  \mathcal{F}_{t}\right)  $-stopping time $T$, and the
optional section theorem.
\end{proof}

\bigskip

In the literature about enlargements of filtrations (\cite{jeulin},
\cite{jeulinyor}, \cite{zurich}, etc.), a number of explicit computations of
supermartingales associated to various $L^{\prime}s$ have been given. We shall
use some of these computations to produce some examples of pseudo-stopping
times, with the help of the proposition.

\begin{enumerate}
\item First let us check again that we recover the example of D. Williams from
the Proposition 5. With the notations of the introduction $\left(
L=\sigma\right)  $, it is not hard to see that (see \cite{rogerswilliams})%
\[
Z_{t}=1-B_{t\wedge T_{1}}^{+}.
\]
Hence
\[
\rho=\sup\left\{  s<\sigma:\text{ }B_{s}=S_{s}\right\}  .
\]

\item Consider $\left(  R_{t}\right)  _{t\geq0}$ a three dimensional Bessel
process, starting from zero, its filtration $\left(  \mathcal{F}_{t}\right)
$, and
\[
L=L_{1}=\sup\left\{  t:\quad R_{t}=1\right\}  .
\]
Then
\begin{equation}
\rho=\sup\left\{  t<L:\quad R_{t}=\sup_{u\leq L}R_{u}\right\}  ,\label{eight}%
\end{equation}
is a $\left(  \mathcal{F}_{t}\right)  $ pseudo-stopping time. This follows
from the fact that
\[
Z_{t}^{L}=1\wedge\frac{1}{R_{t}},
\]
hence (\ref{eight}) is equivalent to:%
\[
\rho=\sup\left\{  t<L:\qquad Z_{t}^{L}=\inf_{u\leq L}Z_{u}^{L}\right\}  ,
\]
and from the above proposition:%
\[
Z_{t}^{\rho}=1\wedge\left(  \frac{1}{\underset{u\leq t}{\sup}R_{u}}\right)  .
\]
We can generalize further this example by noticing that for $n>2$, we have for
$\left(  R_{t}\right)  _{t\geq0}$ a BES$\left(  n\right)  $, $Z_{t}%
^{L}=1\wedge\left(  \frac{1}{R_{t}}\right)  ^{n-2}$.

\item Consider $\left(  B_{u}\right)  _{u\geq0}$\ a one dimensional Brownian
motion, $\left(  \mathcal{F}_{t}\right)  $\ its filtration, and
\[
g_{t}=\sup\left\{  s<t:\quad B_{s}=0\right\}  ,
\]
then
\begin{equation}
\rho_{t}=\sup\left\{  s<g_{t}:\quad\frac{\left\vert B_{s}\right\vert }%
{\sqrt{t-s}}=\sup_{u<g_{t}}\frac{\left\vert B_{u}\right\vert }{\sqrt{t-u}%
}\right\} \label{nine}%
\end{equation}
is a $\mathcal{F}_{t}$\ pseudo-stopping time. Again, this follows from the
fact that $\rho_{t}$\ is in fact defined from $g_{t}\left(  =L\right)  $ as in
the framework preceding the proposition, since:%
\[
Z_{u}^{g_{t}}\equiv\Phi\left(  \frac{\left\vert B_{u}\right\vert }{\sqrt{t-u}%
}\right)  ,
\]
with $\Phi\left(  x\right)  =\mathbb{P}\left(  \left\vert N\right\vert \geq
x\right)  $, where $N$\ is a standard Gaussian.

\item We can reinterpret the previous example via a deterministic time-change.
We remark that we can write:%
\[
\frac{B_{u}}{\sqrt{1-u}}=Y_{\log\frac{1}{1-u}},
\]
where $\left(  Y_{s}\right)  _{s\geq0}$, is an Ornstein-Uhlenbeck process
satisfying
\[
Y_{s}=\beta_{s}+\frac{1}{2}\int_{0}^{s}duY_{u}.
\]
We then deduce from example 3 that
\[
\rho^{^{\prime}}=\sup\left\{  s<L_{0}^{^{\prime}}:\text{\qquad}\left\vert
Y_{s}\right\vert =\sup_{u\leq L_{0}^{^{\prime}}}\left\vert Y_{u}\right\vert
\right\}
\]
is a $\left(  \mathcal{F}_{t}^{^{\prime}}\right)  $ pseudo-stopping time,
where
\[
L_{0}^{^{\prime}}\equiv\log\left(  \frac{1}{1-g_{1}}\right)  =\sup\left\{
s>0,\qquad Y_{s}=0\right\}
\]
and $\left(  \mathcal{F}_{t}^{^{\prime}}\right)  $\ is the natural filtration
of $\left(  Y_{t}\right)  $.

\item Let us consider the case of a transient diffusion $X_{t}$. Let $s$\ be a
scale function such that $s\left(  -\infty\right)  =0$ and $s\left(  x\right)
>0$. Let
\[
L_{a}=\sup\left\{  t;\quad X_{t}=a\right\}  ,
\]
the last passage at the level $a$. We have (see \cite{pitmanyor}):%
\[
Z_{t}^{L_{a}}=1\wedge\frac{s\left(  X_{t}\right)  }{s\left(  a\right)  }.
\]
Thus
\[
\rho_{a}=\sup\left\{  t<L_{a}:\quad s\left(  X_{t}\right)  =\inf_{u\leq L_{a}%
}s\left(  X_{u}\right)  \right\}
\]
is a pseudo-stopping time in the filtration of $\left(  X_{t}\right)  $. For
example, let us consider the case of a brownian motion with a negative drift:
\[
X_{t}\equiv x+\mu t+\sigma B_{t},\quad\mu<0.
\]
In this case, the scale function is
\[
s\left(  x\right)  =\exp\left(  -\frac{2\mu x}{\sigma^{2}}\right)  .
\]
Hence
\[
\rho_{a}=\sup\left\{  t<L_{a}:\quad\mu t+\sigma B_{t}=\inf_{u\leq L_{a}%
}\left(  \mu u+\sigma B_{u}\right)  \right\}
\]
is a pseudo-stopping time in the filtration of $\left(  B_{t}\right)  $.\bigskip
\end{enumerate}

As for D. Williams' example, none of these pseudo-stopping times remains a
pseudo-stopping time in the larger filtration $\left(  \mathcal{F}_{t}%
^{L}\right)  $.This is a consequence of a result of Az\'{e}ma (\cite{azema}).

\begin{proposition}
Let $L$\ be the end of a predictable set such that $\mathbb{P}\left[
L=T\right]  =0$. Then $L$\ is not a pseudo-stopping time.
\end{proposition}

\begin{proof}
From a result of Az\'{e}ma (\cite{azema}), as $A_{t}^{L}=a_{t}^{L}$ is
continuous, the law of $A_{\infty}^{L}$\ is the exponential law of parameter
$1$, whilst for pseudo-stopping times, the law of $A_{\infty}^{L}$\ is
$\delta_{1}$, the Dirac mass at one. Hence $L$\ cannot be a pseudo-stopping time.
\end{proof}

\subsection{Another generalization}

We now give a generalization of the previous construction. We make the same
assumptions about the filtration $\left(  \mathcal{F}_{t}\right)  $\ and the
time $L$, with the extra assumption that $\mathbb{P}\left[  L=\infty\right]
=0$. Let $\left(  \Delta_{t}\right)  $\ be a nonincreasing, continuous and
adapted process such that
\begin{align}
\Delta_{0}  & =1\label{ncond}\\
\Delta_{\infty}  & =0.
\end{align}
Let us define $\rho$\ by%
\[
\rho\equiv\sup\left\{  t<L;\quad Z_{t}=\Delta_{t}\right\}  .
\]
Again, for every $\left(  \mathcal{F}_{t}\right)  $ stopping time $T$, we have%
\[
\left\{  T<\rho\right\}  =\left\{  T^{^{\prime}}<L\right\}
\]
where%
\[
T^{^{\prime}}=\inf\left\{  t>T,\quad Z_{t}\leq\Delta_{T}\right\}
\]
Thus%
\[
\mathbb{E}\left[  Z_{T}^{\rho}\right]  =\mathbb{P}\left[  T<\rho\right]
=\mathbb{P}\left[  T^{^{\prime}}<L\right]  =\mathbb{E}\left[  Z_{T^{^{\prime}%
}}\right]  =\mathbb{E}\left[  \Delta_{T}\right]  ,
\]
and with the optional section theorem we can conclude that
\[
Z_{t}^{\rho}=\Delta_{t},\qquad t\geq0.
\]
It follows from Theorem \ref{mainthm} that $\rho$\ is a pseudo-stopping time.
Hence we have proved the following:

\begin{proposition}
Let $\left(  \Delta_{t}\right)  $\ be a nonincreasing, continuous and adapted
process such that:%
\begin{align*}
\Delta_{0}  & =1\\
\Delta_{\infty}  & =0.
\end{align*}
Then, under the assumptions made above, there always exists a pseudo-stopping
time $\rho$\ such that $Z_{t}^{\rho}=\Delta_{t}$, for $t\geq0$.
\end{proposition}

So we can associate a pseudo-stopping time to any continuous, nonincreasing
adapted process $\left(  \Delta_{t}\right)  $\ which satisfies (\ref{ncond}).
But there is not uniqueness since we can use different $Z^{\prime}s$
associated to different $L^{\prime}s$ to construct $\rho$. In other words,
every continuous, nonincreasing adapted process $\left(  \Delta_{t}\right)  $
satisfying (\ref{ncond}) is the dual predictable projection of some
$1_{\left\{  \rho\leq t\right\}  }$, where $\rho$\ is a pseudo-stopping time.

As an example, we can take
\[
\Delta_{t}=\exp\left(  -S_{t}\right)
\]
with the already used\ notations. Then,
\[
\rho=\sup\left\{  t<\sigma;\quad1-B_{t}^{+}=\exp\left(  -S_{t}\right)
\right\}
\]
is a pseudo-stopping time in the filtration of the Brownian motion $\left(
B_{t}\right)  $. We could as well take
\[
\Delta_{t}=\exp\left(  -L_{t}\right)  ,
\]
where $L_{t}$\ is the Brownian local time at level zero. In that case,
\[
\rho=\sup\left\{  t<\sigma;\quad1-B_{t}^{+}=\exp\left(  -L_{t}\right)
\right\}
\]
is a pseudo-stopping time.

We can also notice that if we take some deterministic $\Delta$, this
construction allows us to construct a pseudo-stopping time with a given
distribution. For example,
\[
\rho=\sup\left\{  t<\sigma;\quad1-B_{t}^{+}=\exp\left(  -\lambda t\right)
\right\}  ,
\]
where $\lambda>0$. Then $\rho$\ follows an exponential law of parameter
$\lambda$.

In the following section, we will see that we can drop the continuity
assumption but we will have to enlarge the initial probability space.

\subsection{Further examples}

In this section, we shall link pseudo-stopping times with other random times
that appear in the literature. In particular, we will see that the random
times allowing the $\left(  \mathbf{H}\right)  $\ hypothesis (see
\cite{elliotjeanbyor}) to hold are special cases of pseudo-stopping times.

\subsubsection{The hypothesis $\left(  \mathbf{H}\right)  $}

First, we give the following obvious result:

\begin{proposition}
If $\rho$\ is a random time that is independent from $\mathcal{F}_{\infty}$,
then it is a pseudo-stopping time.
\end{proposition}

\begin{example}
If $\rho$\ is an exponential time of parameter $\lambda$\ that is independent
from $\mathcal{F}_{\infty}$, then it is a pseudo-stopping time.
\end{example}

\begin{example}
Another example is given by what D. Williams (\cite{williams}) calls a "silly"
time:%
\[
\rho=\frac{1}{1+\left\vert B_{2}-B_{1}\right\vert },
\]
which is independent from $\mathcal{F}_{1}$.\bigskip
\end{example}

Now suppose that our probability space supports a uniform random variable
$\Theta$\ on $\left(  0,1\right)  $ that is independent of the sigma field
$\mathcal{F}_{\infty}$. Assume we are given an $\left(  \mathcal{F}%
_{t}\right)  $-adapted increasing and continuous process satisfying $A_{0}=0$
and $A_{\infty}=1$\ . Let us consider the random time defined by:%
\[
\rho=\inf\left\{  t;\quad A_{t}>\Theta\right\}  .
\]
It is not difficult to check that
\begin{equation}
\mathbb{P}\left[  \rho>t\mid\mathcal{F}_{t}\right]  =1-A_{t}%
.\label{indeppseudo}%
\end{equation}
We have thus constructed a pseudo-stopping time associated with a given
continuous process $\left(  A_{t}\right)  $. This construction is well known,
see \cite{jenbrutk1} for more details and references. But this construction is
not always possible (for example when $\mathcal{F}_{\infty}=\mathcal{F}$),
which explains why our construction in the previous section is more general.

But the pseudo-stopping times that are constructed in the way of
(\ref{indeppseudo}) enjoy the following noticeable property (\cite{jenbrutk1},
\cite{delmeyfil}):
\begin{equation}
\mathbb{P}\left[  \rho>t\mid\mathcal{F}_{t}\right]  =\mathbb{P}\left[
\rho>t\mid\mathcal{F}_{\infty}\right]  .\label{H}%
\end{equation}
Random times with this property are often used in the literature on default
modelling (see \cite{jenbrutk1}, \cite{elliotjeanbyor}) and were studied in
\cite{delmeyfil}, \cite{bremaudyor}. There are several equivalent formulations
for (\ref{H}). Before we mention them, let us notice that any random time
satisfying (\ref{H}) is a pseudo-stopping time. In fact, we have a stronger
result: every $\left(  \mathcal{F}_{t}\right)  $ martingale is an $\left(
\mathcal{F}_{t}^{\rho}\right)  $ martingale (see \cite{delmeyfil}). Thus the
fourth statement in Theorem \ref{mainthm} is satisfied.\bigskip

Now let us consider the $\left(  \mathbf{H}\right)  $ hypothesis in our
framework of progressive enlargement with a random time $\rho$: every $\left(
\mathcal{F}_{t}\right)  $-square integrable martingale is an $\left(
\mathcal{F}_{t}^{\rho}\right)  $-square integrable martingale. This hypothesis
was studied by Dellacherie and Meyer \cite{delmeyfil}, Br\'{e}maud and Yor
\cite{bremaudyor}. It is equivalent to one of the following hypothesis (see
\cite{elliotjeanbyor} for more references):

\begin{enumerate}
\item $\forall t$, the $\sigma$-algebras $\mathcal{F}_{\infty}$\ and
$\mathcal{F}_{t}^{\rho}$\ are conditionally independent given $\mathcal{F}%
_{t}$.

\item For all bounded $\mathcal{F}_{\infty}$-measurable random variables
$\mathbf{F}$\ and all bounded $\mathcal{F}_{t}^{\rho}$-measurable random
variables $\mathbf{G}_{t}$, we have%
\[
\mathbb{E}\left[  \mathbf{FG}_{t}\mid\mathcal{F}_{t}\right]  =\mathbb{E}%
\left[  \mathbf{F}\mid\mathcal{F}_{t}\right]  \mathbb{E}\left[  \mathbf{G}%
_{t}\mid\mathcal{F}_{t}\right]  .
\]

\item For all bounded $\mathcal{F}_{t}^{\rho}$-measurable random variables
$\mathbf{G}_{t}$:%
\[
\mathbb{E}\left[  \mathbf{G}_{t}\mid\mathcal{F}_{\infty}\right]
=\mathbb{E}\left[  \mathbf{G}_{t}\mid\mathcal{F}_{t}\right]  .
\]

\item For all bounded $\mathcal{F}_{\infty}$-measurable random variables
$\mathbf{F}$,
\[
\mathbb{E}\left[  \mathbf{F}\mid\mathcal{F}_{t}^{\rho}\right]  =\mathbb{E}%
\left[  \mathbf{F}\mid\mathcal{F}_{t}\right]  .
\]

\item For all $s\leq t$,
\[
\mathbb{P}\left[  \rho\leq s\mid\mathcal{F}_{t}\right]  =\mathbb{P}\left[
\rho\leq s\mid\mathcal{F}_{\infty}\right]  .
\]
\end{enumerate}

Thus, pseudo-stopping times may be considered as a generalized or a weakened
form of the $\left(  \mathbf{H}\right)  $ hypothesis since then local
martingales in the initial filtration remain local martingales in the enlarged
one up to time $\rho$. Moreover, for most of the examples we have considered,
such as D. Williams', (\ref{H}) is not satisfied.

\subsubsection{Randomized stopping times and F\"{o}llmer measures}

Now we give a relation between pseudo-stopping times and randomized stopping
times as presented in \cite{meyflou}. First we give some definitions. We
always consider a given probability space $\left(  \Omega,\mathcal{F},\left(
\mathcal{F}_{t}\right)  _{t\geq0},\mathbb{P}\right)  $.

\begin{definition}
A randomized random variable on $\left(  \Omega,\mathcal{F},\mathbb{P}\right)
$ is a probability measure $\mu$\ on $\left(  \left[  0,\infty\right]
\times\Omega,\mathcal{B}\left(  \left[  0,\infty\right]  \right)
\otimes\mathcal{F}\right)  $ such that its projection on $\Omega$\ is equal to
$\mathbb{P}$.
\end{definition}

For example, let $\rho$\ be a random time; then $\mu_{\rho}$ defined by
\[
\mu_{\rho}\left(  X\right)  =\mathbb{E}\left[  X_{\rho}\right]  ,
\]
for all bounded measurable processes $\left(  X_{t}\right)  $ is a randomized
random variable. \bigskip

We know from a result of F\"{o}llmer (see \cite{dellachmeyer}) that there
exists an increasing c\`{a}dl\`{a}g process $\left(  A_{t}\right)  $ such that
$A_{0}=0$ and
\[
\mu\left(  X\right)  =\mathbb{E}\left[  \int_{0}^{\infty}X_{s}dA_{s}\right]  ,
\]
for all nonnegative process $\left(  X_{t}\right)  $. The fact that the
projection on $\Omega$\ is equal to $\mathbb{P}$ means that $A_{\infty}=1,$
$a.s$.

\begin{definition}
If the process $\left(  A_{t}\right)  $ associated with $\mu$\ on $\left(
\left[  0,\infty\right]  \times\Omega,\mathcal{B}\left(  \left[
0,\infty\right]  \right)  \otimes\mathcal{F}\right)  $ is adapted, then we say
that $\mu$\ is a randomized stopping time.
\end{definition}

By considering the new space $\overline{\Omega}=\left[  0,1\right]
\times\Omega$ endowed with the $\sigma$-fields $\overline{\mathcal{F}}%
$=$\mathcal{B}\left(  \left[  0,1\right]  \right)  \otimes\mathcal{F}$,
$\overline{\mathcal{F}}_{t}$=$\mathcal{B}\left(  \left[  0,1\right]  \right)
\otimes\mathcal{F}_{t}$ (augmented in the usual way) and the probability
measure $\overline{\mathbb{P}}=\lambda\otimes\mathbb{P}$, it is possible to
show that for every randomized stopping time $\mu$, there exists a stopping
time $\rho$\ in this new filtered space such that
\[
\mu\left(  X\right)  =\overline{\mathbb{E}}\left[  X_{\rho}\right]  ,
\]
for all bounded measurable process $\left(  X_{t}\right)  $ on $\left(
\left[  0,\infty\right]  \times\Omega,\mathcal{B}\left(  \left[
0,\infty\right]  \right)  \otimes\mathcal{F}\right)  $. We take the convention
that a random variable $H$\ on $\Omega$ can be considered as the random
variable on $\overline{\Omega}$: $\left(  u,\omega\right)  \rightarrow
H\left(  \omega\right)  $. Conversely to every stopping time of $\overline
{\mathcal{F}}_{t}$ corresponds a randomized stopping time. \bigskip

This construction is always carried on the enlarged space $\overline{\Omega}
$. The third statement in Theorem \ref{mainthm} allows us to use
pseudo-stopping times to construct randomized stopping times without enlarging
the initial space.

\begin{proposition}
Let $\rho$\ be a pseudo-stopping time and $A_{t}^{\rho}$\ the $\left(
\mathcal{F}_{t}\right)  $\ dual optional projection of the process
$1_{\left\{  \rho\leq t\right\}  }$. Then the F\"{o}ellmer measure $\mu
$\ associated with $A_{t}^{\rho}$ is a randomized stopping time. Moreover, for
every bounded or nonnegative $\left(  \mathcal{F}_{t}\right)  $\ optional
process $\left(  X_{t}\right)  $:%
\[
\mu\left(  X\right)  =\mathbb{E}\left[  X_{\rho}\right]  .
\]
\end{proposition}

\subsubsection{Randomized stopping times and families of stopping times}

\begin{proposition}
Let $\left(  T_{u}\right)  _{u\geq0}$ be a family of $\left(  \mathcal{F}%
_{t}\right)  $\ stopping times and $S$ a positive random variable, independent
of the family $\left(  \mathcal{F}_{\infty}\right)  $. Then%
\[
\rho=T_{S}%
\]
is a $\left(  \mathcal{F}_{t}\right)  $\ pseudo-stopping time.
\end{proposition}

\begin{proof}
Let $\left(  M_{t}\right)  $ be a bounded $\left(  \mathcal{F}_{t}\right)
$\ martingale;
\begin{align*}
\mathbb{E}\left[  M_{T_{S}}\right]   & =\mathbb{E}\left[  \mathbb{E}\left[
M_{T_{s}}\mid S=s\right]  \right] \\
& =\mathbb{E}\left[  \mathbb{E}\left[  M_{0}\right]  \mid S=s\right] \\
& =\mathbb{E}\left[  M_{0}\right]  .
\end{align*}
\end{proof}

The previous proposition shows that any independently time changed family of
stopping times is a pseudo-stopping time. In fact, this proposition admits a
converse: every pseudo-stopping time is, in law, a time changed family of
stopping times. More precisely:

\begin{proposition}
Let $\rho$\ be a $\left(  \mathcal{F}_{t}\right)  $\ pseudo-stopping time,
which avoids all $\left(  \mathcal{F}_{t}\right)  $-stopping times, and
$Z_{t}=\mathbb{P}\left[  \rho>t\mid\mathcal{F}_{t}\right]  $ its associated
supermartingale. Set
\[
\alpha_{u}\equiv\inf\left\{  t\geq0,\quad\left(  1-Z_{t}\right)  >u,\quad0\leq
u\leq1\right\}  ,
\]
the right-continuous generalized inverse of the increasing continuous process
$\left(  1-Z_{t}\right)  $. Then $\left(  \alpha_{u}\right)  _{0\leq u\leq1}%
$\ is a family of $\left(  \mathcal{F}_{t}\right)  $\ stopping times and%
\[
\rho\overset{law}{=}\alpha_{U},
\]
where $U$\ is a random variable with uniform law, independent of $\left(
\mathcal{F}_{\infty}\right)  $.
\end{proposition}

\begin{proof}
The fact $\alpha_{u}$ is a stopping time, for all $u$, follows from
\[
\left\{  \alpha_{u}\leq t\right\}  =\left\{  u\leq\left(  1-Z_{t}\right)
\right\}  ,\quad\forall t\geq0.
\]
From (\ref{formulerandom}), we also have
\[
\mathbb{E}\left[  g\left(  \rho\right)  \right]  =\mathbb{E}\left[  \int
_{0}^{1}g\left(  \alpha_{u}\right)  du\right]  ,
\]
for all bounded Borel function $g$. This establishes $\left(  \rho
\overset{law}{=}\alpha_{U}\right)  $.
\end{proof}

\section{A discrete analogue: the coin-tossing case}

Let $\left(  X_{n}\right)  _{n\geq1}$ be the standard random walk with
Bernoulli increments. In his paper \cite{Le Gall}, Le Gall proved an analogue
of Williams' path decomposition for $\left(  X_{n}\right)  $. To fix ideas, we
shall consider the canonical space $\Omega=\mathbb{Z}^{N}$ endowed with the
product $\sigma$-field. $\left(  X_{n}\right)  $ will be the coordinate
process and $\left(  \mathbb{P}_{x}\right)  _{x\in\mathbb{Z}}$\ the family of
probability laws which make $\left(  X_{n}\right)  $\ the standard random walk
with Bernoulli increments. We also denote by $\left(  \mathbb{Q}_{x}\right)
_{x\in N}$\ the unique family of probability measures such that $\left(
X_{n},\mathbb{Q}_{x}\right)  $\ is a Markov chain with transition
probabilities:%
\begin{align*}
\mathbb{Q}_{0}\left[  X_{1}=1\right]   & =1\\
if\text{ }x  & \geq1,\quad\mathbb{Q}_{x}\left[  X_{1}=x+1\right]  =\frac{1}%
{2}\left(  1+\frac{1}{x}\right)  ;\quad\mathbb{Q}_{x}\left[  X_{1}=x-1\right]
=\frac{1}{2}\left(  1-\frac{1}{x}\right)  .
\end{align*}
Now let $p\geq1$ and define:%
\begin{align*}
\sigma_{p}  & =\inf\left\{  k;\quad X_{k}=p\right\}  ,\\
\eta & =\sup\left\{  k\leq\sigma_{p}:\quad X_{k}=0\right\}  ,\\
m  & =\sup\left\{  X_{k},\quad k\leq\eta\right\}  ,\\
\gamma & =\inf\left\{  k\geq0;\quad X_{k}=m\right\}  .
\end{align*}
Then, Le Gall's statement is that under $\mathbb{P}_{0}$:

\begin{enumerate}
\item The processes $\left(  X_{k}\right)  _{0\leq k\leq\eta}$ and $\left(
X_{\eta+k}\right)  _{0\leq k\leq\sigma_{p}-\eta}$ are independent, with the
second being distributed as $\left(  X_{k}\right)  _{0\leq k\leq\sigma_{p}}$
under $\mathbb{Q}_{0}$;

\item $m$ is uniformly distributed on $\left\{  0,1,\ldots,p-1\right\}  $;

\item Conditionally on $\left\{  m=j\right\}  $, the processes $\left(
X_{k}\right)  _{0\leq k\leq\gamma}$ and $\left(  X_{\eta-k}\right)  _{0\leq
k\leq\eta-\gamma}$ are independent, the first being distributed as $\left(
X_{k}\right)  _{0\leq k\leq\sigma_{j}}$ under $\mathbb{P}_{0}$, and the second
as $\left(  X_{k}\right)  _{0\leq k\leq\sigma_{j+1}-1}$ under $\mathbb{Q}_{0}$.
\end{enumerate}

\begin{proposition}
If $\left(  M_{n}\right)  _{n\in\mathbb{N}}$ is a bounded martingale, then
\[
\mathbb{E}_{0}\left[  M_{\gamma}\right]  =\mathbb{E}_{0}\left[  M_{0}\right]
.
\]
Thus $\gamma$\ is a pseudo-stopping time.
\end{proposition}

\begin{proof}
The discrete time setup allows us to give an elementary argument, based in
part on the fact that $M_{n}$, as every $\mathcal{F}_{n}$ measurable variable,
may be written as:%
\[
M_{n}=f_{n}\left(  X_{1},X_{2},\ldots,X_{n}\right)  ,
\]
where $f_{n}$ is a bounded function depending on $n$\ variables.

Now, for any bounded function $g$:%
\[
\mathbb{E}_{0}\left[  M_{\gamma}g\left(  m\right)  \right]  =\mathbb{E}%
_{0}\left[  \mathbb{E}_{0}\left[  M_{\gamma}\mid m\right]  g\left(  m\right)
\right]  .
\]
But, from $\left(  3\right)  $\ in Le Gall's satatement:%
\begin{align*}
\mathbb{E}_{0}\left[  M_{\gamma}\mid m=j\right]   & =\mathbb{E}_{0}\left[
f_{\sigma_{j}}\left(  X_{1},X_{2},\ldots,X_{\sigma_{j}}\right)  \right] \\
& =\mathbb{E}_{0}\left[  M_{\sigma_{j}}\right]  =\mathbb{E}_{0}\left[
M_{0}\right]  .
\end{align*}
Thus, we have obtained:%
\begin{align*}
\mathbb{E}_{0}\left[  M_{\gamma}g\left(  m\right)  \right]   & =\mathbb{E}%
_{0}\left[  M_{\gamma}\right]  \mathbb{E}_{0}\left[  g\left(  m\right)
\right] \\
& =\mathbb{E}_{0}\left[  M_{\infty}\right]  \mathbb{E}_{0}\left[  g\left(
m\right)  \right]  ,
\end{align*}
which is the discrete analogue of Proposition 2, and shows a fortiori that
$\gamma$\ is a pseudo-stopping time.
\end{proof}


\newpage

\begin{thebibliography}{99}                                                                                                %
\bibitem {azema}\textsc{Az\'{e}ma, J.} (1972). Quelques
applications de la th\'{e}orie g\'{e}n\'{e}rale des processus I.
\textit{Invent. Math.} \textbf{18}
 293-336.

\bibitem {barlyor}\textsc{Barlow, M.T. and Yor, M.} (1981). (Semi)-martingale
inequalities and local times. \textit{Z.f.W.} \textbf{55} 237-254.

\bibitem {bremaudyor}\textsc{Br\'{e}maud, P. and Yor, M.} (1978). Changes of
filtration and of probability measures. \textit{Z.f.W.} \textbf{45}
269-295.

\bibitem {delmaismey}\textsc{Dellacherie, C., Maisonneuve, B. and
Meyer, P.A.} (1992). \textit{Probabilit\'{e}s et potentiel.
Chapitres XVII-XXIV: Processus de Markov (fin), Compl\'{e}ments de
calcul stochastique.}
 Hermann.

\bibitem {delmeyfil}\textsc{Dellacherie, C. and Meyer, P.A.} (1978). A propos du
travail de Yor sur les grossissements des tribus.
\textit{S\'{e}m.Proba. XII, Lecture Notes in Mathematics.}
\textbf{649} 69-78.


\bibitem {dellachmeyer}\textsc{Dellacherie, C. and Meyer, P.A.}:
\textit{Probabilit\'{e}s et potentiel}, Hermann, Paris, vol. I 1976,
vol. II 1980.

\bibitem {elliotjeanbyor}\textsc{R.J. Elliott, R.J., Jeanblanc, M. and Yor,
M.}(2000). On models of default risk. \textit{Math. Finance.}
\textbf{10} 179-196.

\bibitem {jenbrutk1}\textsc{Jeanblanc, M. and Rutkowski, M.} (2000). Modeling
default risk: Mathematical tools. Fixed Income and Credit risk
modeling and Management, New York University, Stern school of
business, Statistics and operations research department, Workshop .

\bibitem {jeulin}\textsc{Jeulin, T.} (1980). \textit{Semi-martingales et
grossissements d'une filtration}, Lecture Notes in Mathematics
\textbf{833}, Springer.

\bibitem {yorjeulin}\textsc{Jeulin, T. and Yor, M.} (1978a). Grossissement d'une
filtration et semimartingales: formules explicites.
\textit{S\'{e}m.Proba. XII, Lecture Notes in Mathematics}
\textbf{649} 78-97.

\bibitem {jeulinyor}\textsc{Jeulin, T. and Yor, M. (eds).} (1985b). \textit{Grossissements
de filtrations: exemples et applications}, Lecture Notes in
Mathematics \textbf{1118}, Springer.

\bibitem {Knight-Maisonneuve}\textsc{Knight, F.B. and Maisonneuve, B.}(1994). A
characterization of stopping times. \textit{Annals of probability.}
\textbf{22} 1600-1606.

\bibitem {Le Gall}\textsc{Le Gall, J.F.} (1984). Une approche
\'{e}l\'{e}mentaire des th\'{e}or\`{e}mes de d\'{e}composition \ de
Williams. \textit{S\'{e}m.Proba. XX, Lecture Notes in Mathematics.}
\textbf{1204} 447-464.

\bibitem {meyflou}\textsc{Meyer, P.A.} (1978). Convergence faible et
compacit\'{e} des temps d'arr\^{e}t, d'apr\`{e}s
Baxter-Chacon.\textit{ S\'{e}m.Proba. XII, Lecture Notes in
Mathematics.} \textbf{649} 411-424.

\bibitem {pitmanyor}\textsc{Pitman, J.W and Yor, M.} (1981). Bessel processes and
infinitely divisible laws. In: D. Williams (ed.) Stochastic
integrals, Lecture Notes in Mathematics \textbf{851}, Springer.

\bibitem {protter}\textsc{Protter, P.E.} (2003). \textit{Stochastic integration and
differential equations.} Springer, Second edition.

\bibitem {revuzyor}\textsc{Revuz, D. and Yor, M.}(1999). \textit{Continuous martingales
and Brownian motion.} Springer, Third edition .

\bibitem {rogerswilliams}\textsc{Rogers,C. and Williams, D.}(1987). \textit{Diffusions,
Markov processes and Martingales, vol 2: Ito calculus.} Wiley and
Sons, New York.

\bibitem {williams}\textsc{ Williams, D.} (2002). A non stopping time with the
optional-stopping property. \textit{Bull. London Math. Soc.}
\textbf{34}, 610-612.

\bibitem {zurich}\textsc{Yor, M.} (1997). \textit{Some aspects of Brownian motion,
Part II. Some recent martingales problems.} Birkhauser, Basel.
\end{thebibliography}
\end{document}